\documentclass[12pt]{article}

\usepackage[utf8]{inputenc}

\usepackage{amsbsy, amsmath,amsfonts,amssymb,amsthm,amscd,amssymb,latexsym,float}
\usepackage{graphicx}
\usepackage[dvips]{epsfig}
\usepackage{graphics}
\usepackage{psfrag}

\usepackage[usenames,dvipsnames]{xcolor}

 \usepackage{hyperref}

 \theoremstyle{plain}

\newtheorem{remark}{Remark}
\theoremstyle{definition}
\newtheorem{definition}{Definition}
\newtheorem{theorem}{Theorem}

\newtheorem{conclusion}[theorem]{Conclusion}

\newtheorem{conjecture}[theorem]{Conjecture}

\newtheorem{problem}{Problem}

\newcommand{\R}{\mathbb{R}}

 \providecommand{\keywords}[1]
 {
 	\small	
 	\textbf{\textit{Keywords---}} #1
 }

\newenvironment{dedication}
{\vspace{6ex}\begin{quotation}\begin{center}\begin{em}}
			{\par\end{em}\end{center}\end{quotation}}
 
\begin{document}
 	
\title{ An encounter 
 in the realm of Structural Stability 
 between  a qualitative theory for geometric shapes and   
 one for the 
 integral foliations  of 
  differential equations 
 }

\author{Jorge Sotomayor\thanks{The  
		author is a 
		fellow of  CNPq,   Grant: PQ-SR- 307690/2016-4.
}}

\maketitle


\begin{dedication}
 To the memory of  Maur\'icio M. Peixoto (1921 - 2019), Carlos E. Harle (1937 - 2020),   Carlos T. Guti\'errez  (1944 - 2008) and Daniel B. Henry (1945 - 2002)
\end{dedication}


  \begin{abstract}
This evocative  essay  focuses on some 
 landmarks that  led the author 
 to the study
of  principal curvature 
 configurations on surfaces in  $\mathbb R^3$, their structural stability
and   generic properties. 
The  starting point was an encounter with  the book of D.  Struik and the reading 
of the  references to the works of  Euler, Monge and Darboux  found  there.
 The  concatenation of these references  with the work of Peixoto, 1962,  on  differential equations  on surfaces, was a crucial second step.
 The  circumstances of the  convergence  toward the theorems of Guti\'errez and Sotomayor, 1982 - 1983, are recounted here. 
  The above 1982 - 1983 theorems are pointed out as  the   first encounter between  the line of thought disclosed 
   from the works of   Monge, 1796, Dupin, 1815,  and Darboux,  1896,   with that transpiring from the achievements  of Poincar\'e, 1881, Andronov - Pontrjagin, 1937, and Peixoto, 1962.
   Some mathematical developments sprouting from the  1982 - 1983  works are mentioned on the final section of this essay.   

   \end{abstract}


 
\maketitle

\section{Monge's Ellipsoid. }\label{sc:1}

This story begins on a hot October night in 1970
in Rio de Janeiro. Suffering  from  insomnia, I decided to snoop around
the books that my wife had carefully accommodated in our 
shelf.  

My candid and relaxed  attitude contrasted with a strange tension
that emanated from the ensemble of books,   
flooding the living  room. Intrigued,  I 
was 
impelled to find out the cause. 

Squeezed  in a corner,
wrapped in an elegant green cover, nervously pulsed the Aguilar, Spanish translation of the  book
``{\it Lectures on Classical Differential Geometry}"\,    of  D. Struik.

\begin{figure}[H]
\begin{center}
\includegraphics[scale=1]{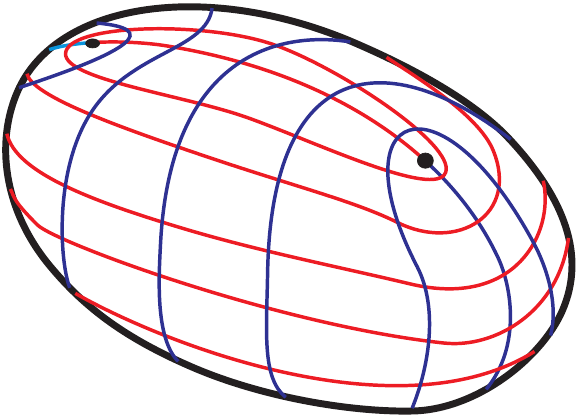}
\end{center}
\caption{\small { Monge's  Ellipsoid.
 Illustration of  the ellipsoid with  three different axes:  
$x^2/a^2  + y^2/b^2 + z^2/c^2 = 1, \,  a> b > c> 0,$
 endowed with its Principal Curvature Configuration, which consists on its umbilic 
singular points and, outside them, its  principal  curvature lines.
  The separatrices  are the curvature lines that connect the umbilic points.} \label{fig:1}}
\end{figure}

My  
devotion to Geometry  and  to  the language of Cervantes, 
drove me to, naively, open that book.  And I did it  just on 
a  page from which, as if it had been lurking,
the
picture of the triaxial ellipsoid, posted in Fig. \ref{fig:1},  popped up.

 More than a century and a half had elapsed
since the 
French mathematician Gaspard Monge conceived it,
calculating its principal  curvature lines, and locating its
 four umbilic points. 
 Meanwhile, printed on books, restricted to an asphyxiating 
two-dimensional existence, the Ellipsoid had traveled
thousands of kilometers, transposed mountains and crossed wide seas,
 so that on that singular  tropical night,  we could meet  face to face.
 
 \begin{figure}[H]
\begin{center}
\includegraphics[scale=.9]{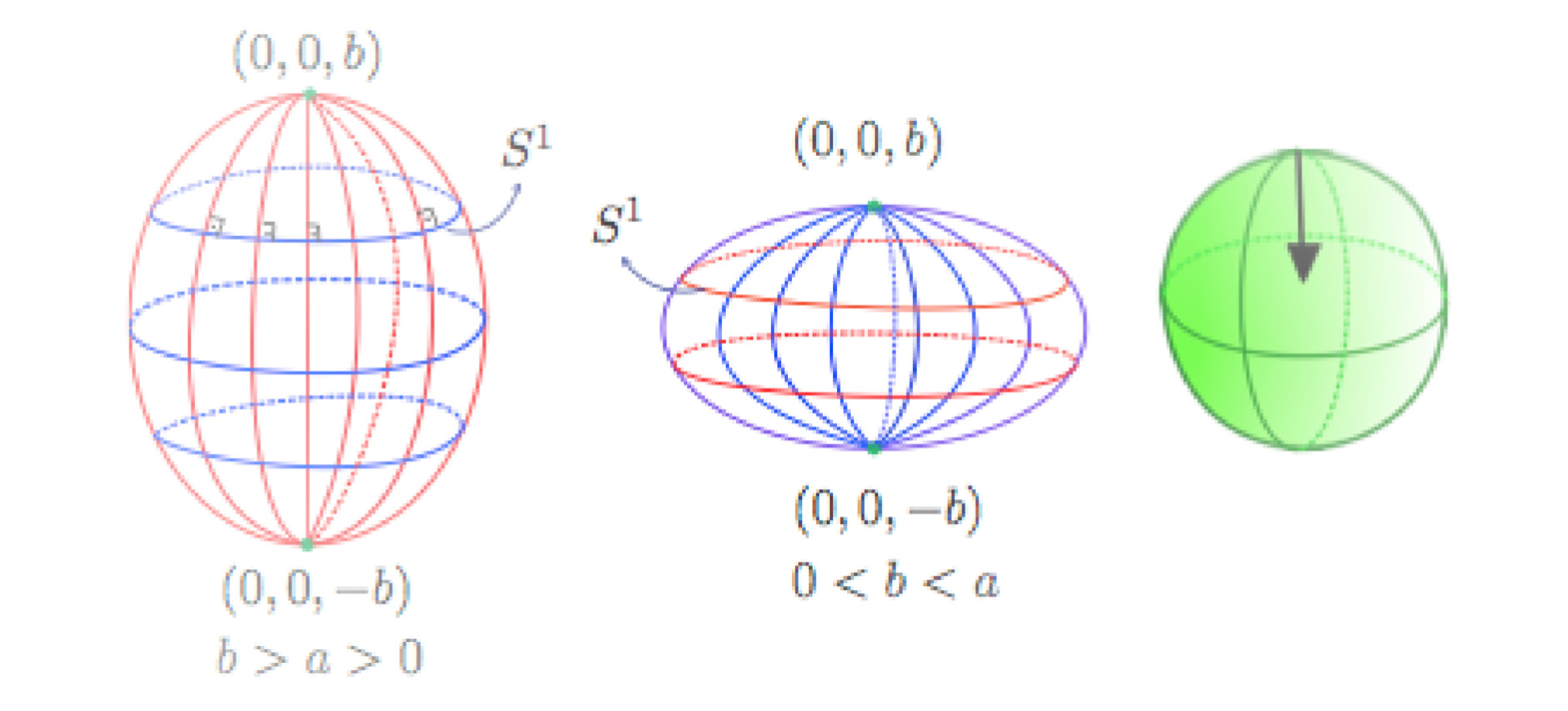}
\end{center}
\caption{{\small Ellipsoid of   Revolution. Illustration of  the (oblate, left) ellipsoid with two equal  axes:  
 $x^2/a^2  + y^2/c^2 + z^2/b^2 = 1, \, 0 < a = c  < b,$
its principal  curvature lines and  its pair of umbilic 
singular points. The (prolate, center) ellipsoid corresponds to $ 0<  b < a =  c.$  On the sphere, right,   $ 0<  c = a =  b, $   all the points are
 umbilic.}  \label{fig:2}}
\end{figure}

As early as in 1961 I had browsed 
through 
Struik's book, but I did not notice the picture. 
 At that time, 
I also was fond to peruse  the classic  ``{\em Geometry and the Imagination}"   of
Hilbert and Cohn Vossen, not so much to read it but  for  contemplating 
its  pictures. 
From this book, as Struik annotates, was taken the
alluded picture. 
 
How could I have   
ignored it?

 With meticulous attention   I  had studied  two other  excellent textbooks, of a comparable level
 to that of Struik:   Willmore's  ``{\it Introduction to
Differential Geometry}" \, 
  as an undergraduate 
  in 1961, and  O'Neill's, ``{\it Elementary
Differential Geometry}",   
as a young instructor in 1967. 
However, none  of these two remarkable models of  geometric  exposition 
contains not  even an outline   of the fascinating picture. 
 It should be mentioned, however,  that
O'Neil proposes a commented exercise in which he outlines  how to locate
the umbilic points of the ellipsoid, with no concern with their local principal configurations. 

\vspace{0.3cm}
 It was like a love at  first glance. The symmetry and beauty of its curves captivated me
immediately. 
As in a ritual of
mutual measuring,
we
stared at each other  for some minutes. 

I promptly read the adjacent text, as well as other
relevant sections.

\vspace{0.3cm}

The reader of that singular night, however,   was somewhat distant from the  
naive 
student of 1961. 
Having traveled
along  pathways of the Qualitative Theory of Differential Equations and Dynamical Systems,
he  had his  
mathematical awareness 
enlarged and shaped by the  study of Structural Stability and
Bifurcation Theory.  
\vspace{0.3cm}

In an hour of active reading I made a mathematical journey that,
chronologically, had spanned  almost two centuries.

There,   the classical  results concerning the  principal curvatures  on a surface  were  proved.  
 I reviewed  Euler's formula that expresses the
normal curvature in terms of the  principal curvatures and directions.

 \begin{figure}[H]
 	\begin{center}
		\includegraphics[scale=1]{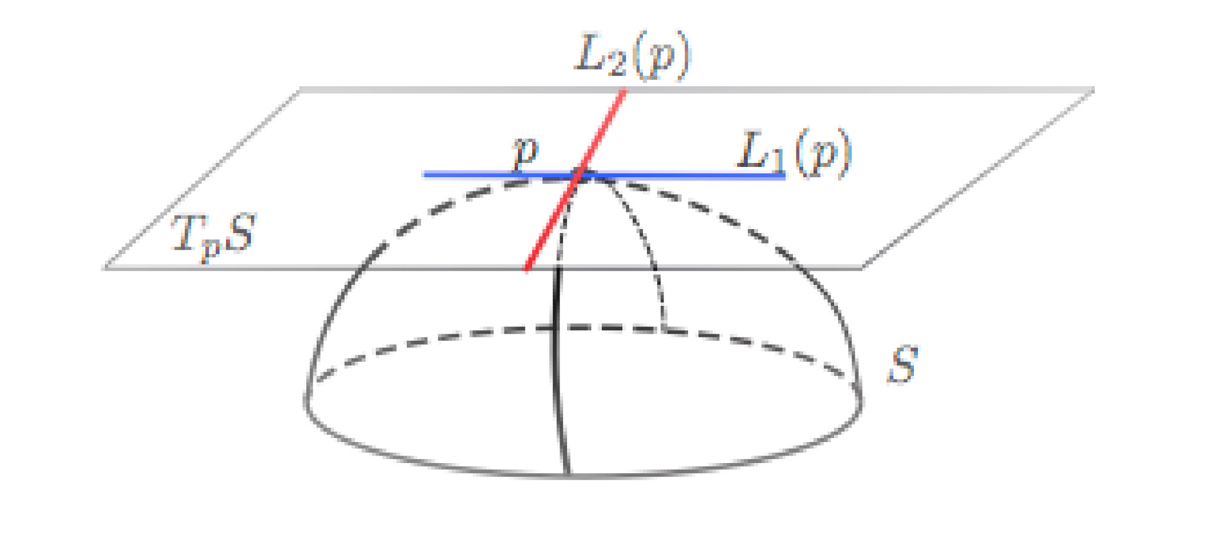}
 	\caption{  
	Principal Directions  at a point $p$ on a Surface $S$ oriented by the inward normal  field. \label{fig:3}}
 	\end{center}
 \end{figure}

The \emph{ principal  curvature lines} are the curves,  along a surface,
 whose tangent directions, \emph{ denominated   principal directions}, make   it to  bend extremely  in  $\mathbb R^ 3 $. The scalar measures  of  these  extreme curvatures, are
 called the \emph{principal curvatures}. 
 One of them, the minimum, is denoted by 
$ k_1 $; the other, the maximum, is designated by $ k_2 $. 
Their  values are given by the normal curvature (the second fundamental form of the
surface),  evaluated at  the  principal directions. So, $ k_1 <
k_2 $, except  at  the \emph{umbilic points}, where $ k_1 = k_2 $.

Euler's formula states that the normal curvature 
$ k_n (\theta) $ in the direction that makes angle $ \theta $ with the  
minimum principal  direction, is given by

$$k_n (\theta) = k_1  \cos^2 \theta + k_2 \sin^2 \theta.$$

This formula is equivalent to the 
 diagonalization,  by means of a rotation, 
of the 
Second  Fundamental   Form of Surface Theory.

\vspace{0.3cm}
In {\it Recherches sur la courbure des surfaces,}  M\'emoires
de  l'Aca\-de\-mie de Sciences de Berlin, 16, (1760),  Euler 
accomplished  
the first study
of these mathematical objects, with which he inaugurated the use of
Differential Calculus in the investigation of surface  geometry. It is derived from this work that, in general,
the   principal  directions are orthogonal.

At  umbilic  points, the   principal  curvatures coincide.
Only outside them, the  principal directions are  
defined and determine  a pair  of 
mutually orthogonal tangent line fields,  $ L_1 $  and  $ L_2$, called the principal line fields;  $ L_1 $, corresponds to the minimal 
and $ L_2 $, to the maximal principal curvatures. The surface is assumed to be oriented. An exchange in its orientation permutes
 $ L_1 $ and  $ L_2$. 

\vspace{0.3cm}
With Struik I reviewed Dupin's Theorem that allows  to determine the
lines of principal curvature of surfaces that belong to
a triply  orthogonal family. 
To this end  it is enough to 
intersect the surface under study with those of the two families
orthogonal to it. 
In the case of Monge's Ellipsoid,  one  must
take  the  hyperboloids of one and two sheets, which, together with
the ellipsoids,   form the so-called   homofocal   system of quadrics. These
quadrics are also the coordinate surfaces of the so-called
Ellipsoidal coordinate system in $ \R^3 $.

\begin{figure}[h]
\begin{center}
\includegraphics[scale=.7]{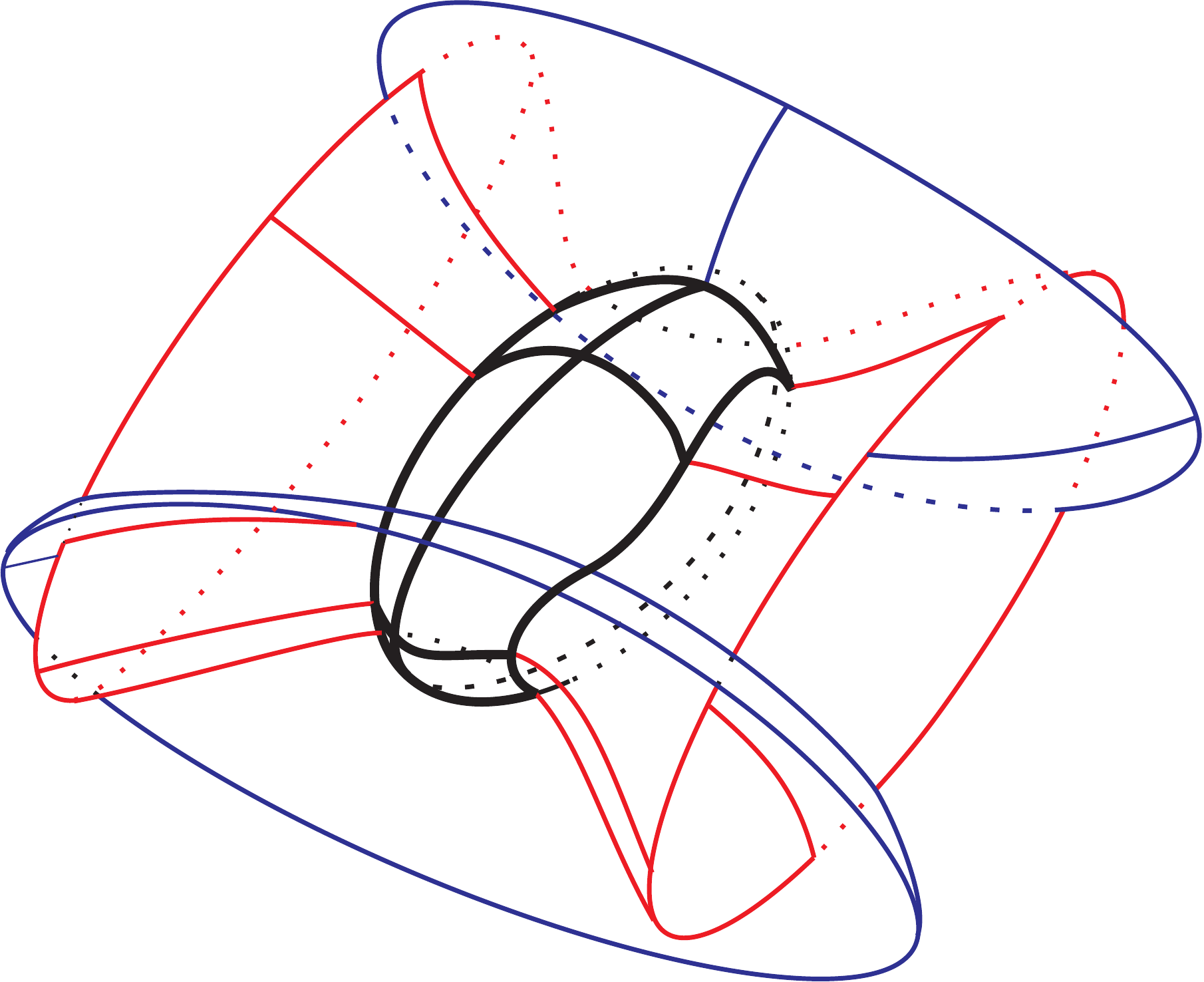}
\end{center}
\caption{ Dupin's Theorem. Illustration  of Principal Configurations for Homofocal Quadrics: 
Ellipsoids and Hyperboloids,  with one and two sheets.   \label{fig:4}}
\end{figure}

\vspace{0.5cm}

\subsection{Report  of a first encounter between Classical Principal Curvature  Geometry and the Structural Stability of Differential Equations.} \label{ss:oencounter_oneiric}
I felt enraptured and 
rewarded.  The browsing  and reading session had allowed me
to become aware of 
a remarkable qualitative jump in the evolution of
the mathematical ideas, in three steps.

 \begin{remark} \label{rm:3jump} 
{\bf A triple qualitative jump. \label{rm:3jump} 
}\\
Firstly, Euler established the definitions
and basic concepts. \\
The second step was taken by Monge,
calculating   a  key example. \\
Thirdly, Dupin demonstrated a general  theorem
that unifies the particular cases, all integrable,  known at that time.
 \end{remark}

From Dupin's Theorem, 
a handy set of color pencils and  numerous diagrams  allowed me to 
 reach 
the following conclusions:  
 
\begin{conclusion}\label{con:oneiric.0}
{\it  Consider the nine dimensional space $ \mathcal Q $ of quadratic compact surfaces, 
 endowed with the topology defined by the, normalized,  coefficients.
The surfaces whose principal  configurations are
Structurally Stable (that is, they are not altered
topologically by small perturbations of their coefficients),
are precisely the ellipsoids of  Monge,  whose three axes are distinct. This class, denoted by $ \mathcal E_3 $, 
is an open and dense  subset in  $\mathcal Q $.}
\end{conclusion}

\begin{conclusion}\label{con:oneiric.2}
{\it 
Inside $\mathcal Q_1 $ = $ \mathcal Q \setminus \mathcal E _3 $,
those whose principal  configurations are First Order Structurally
Stable  (that is, they are not altered by small
disturbances within $ \mathcal Q _1 $), are the ellipsoids of
non-spherical revolution; these form a sub-variety $ \mathcal E
_2 $, of codimension 2 in $\mathcal Q $.}
\end{conclusion}

\begin{conclusion}\label{con:oneiric.6}
{\it  Inside the  complement of $\mathcal E_3 \cup \mathcal E_2$ in  $ \mathcal Q $, the spheres  form a submanifold  of dimension $3$, that is of codimension $6$.}
\end{conclusion}
Analogous results could  also be formulated for the
 non compact quadrics. 
While I was quitting the  long  working session, the following inquiry struck me:

\begin{problem} \label{prob:1}
{\it How would  it be,   in general,  
the principally  structurally
stable  configurations  of  smooth  oriented compact non quadratic  surfaces?}
\end{problem}
 \section{ The Fundamental Problem of  Principal Curvature Geometry.} \label{sc:2}
Two references from  Struik  had attracted my attention. One of them 
was  a note in Volume III of the famous classical  treatise of  Gaston Darboux
``{\it Le\c cons de la The\'orie des Surfaces},"  Gauthier - Villars, 1888 - 1896,   the other was  an
article in Acta Matematica, 1904, by Alvar Gullstrand. 
Having slept over Problem \ref{prob:1}, 
 next  morning,   I  promptly perused   them at the library of IMPA\footnote{National Institute of Pure and Applied Mathematics}. 
  
Both were concerned with the possible configurations of the lines of
principal curvature in the vicinity of an umbilic point.

\begin{figure}[H]
\centerline{
\includegraphics[scale=1.2]{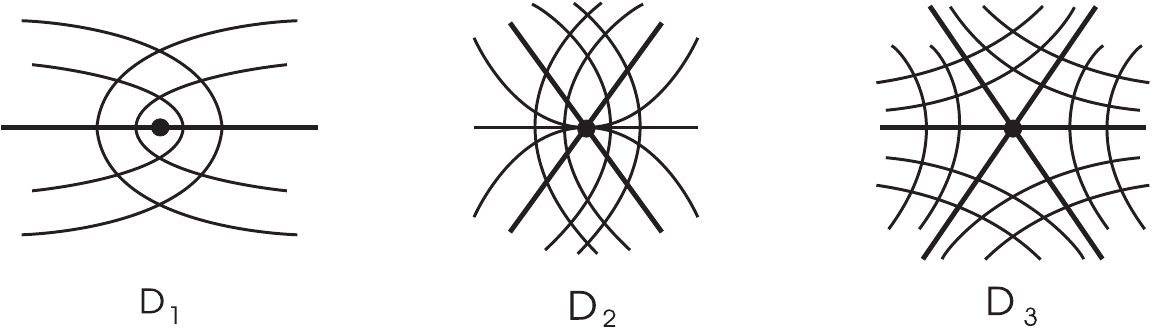}}
\caption{{\small  Darbouxian  Umbilic Points. The subscripts stand for the number of umbilic separatrices 
for each of the principal foliations}.} \label{fig:5}
\end{figure}

The first one had  
a description of the
three cases of generic umbilic points, characterized by
algebraic conditions, expressed as inequalities,  involving   the third derivatives of  
 a function  representing the surface in Monge coordinates
centered at  the umbilic point. See illustrations  in Fig. \ref{fig:5}.
Everything suggested that
principal configurations of these types could be also 
structurally stable, locally at the umbilic point. 
 
 As was  usual, in works  of  around  1886, they were real analytic objects; surfaces in our
case.

The details of the proofs  at the first readings  seemed very  hard to assimilate.

 Consider 
a surface in  Monge form:
$$z = (k/2) (x^2+y^2) + (a/6)x^3 + (b/2)xy^2 + (c/6)y^3 +
O[(x^2+y^2)^2],$$
 \noindent in which the coefficient of the term
$x^2 y$ has been eliminated by means of a rotation.

The conditions that  define the  Darbouxian  umbilic points are written as follows:

\begin{definition} \label{def:darboux}
\vspace{0.4cm}
\item [T)] Transversality Condition:
$b(b-a) \neq 0$; 
\item [D)]    Discriminant Conditions:\newline
D$_{1}$: $
\frac{a}{b}>\Big( \frac{c}{2b} \Big)^{2}+2$; \newline
 D$_{2}$: $\Big( \frac{c}{2b} \Big)^{2}+2> \frac{a}{b}> 1$; $a \neq 2b$; \newline
D$_{3}$: $ \frac{a}{b}<1$.

\end{definition}

\vspace{0.2cm}

These expressions first  appeared  in print in 1983 in  the first step carried out by   Guti\'errez and  Sotomayor, 
cited at the end of 
 section \ref{sc:1981},  in an endevour to solve Problem \ref{prob:1}.  

\vspace{0.2cm}
Being illiterate in German, from the pictures and formulas, I concluded  that  Gullstrand's  work,
after reviewing Darboux's contribution,  proposed a
description of the initial non-generic cases.
 
That seemed  a sort of  the  starting for  an
investigation of the bifurcations of the umbilic points. 
 The potential connections 
with  bifurcation problems of differential equations captured my interest.

 After the contact with  Darboux and  Gullstrand references,  my geometric imagery 
enriched 
remarkably.

 Back home, for a long  while, 
 I contemplated 
 and 
 tried to organize
coherently the 
precious 
antique 
pieces 
I had already collected.

\vspace{0.2cm}
 I  quitted  to rest earlier than usual. 

\vspace{0.2cm}
  
I woke up at dawn and, 
certain of  
being  contributing to the fourth stage  in 
 the   multiple
 mathematical qualitative jump 
 witnessed 
 the  previous  night,  
 annotated  in remark  \ref{rm:3jump},  without hesitation, I  wrote:  
\begin{problem} \label{prob:2} { \it
The fundamental problem on  the study  of  the Principal  Curvature Configurations
of the differentiable, compact and oriented surfaces immersed into $\R^3$, consists
in establishing 
the following two theorems:}
\end{problem}
\begin{theorem} \label{th.1}
The 
necessary 
and 
sufficient 
conditions for a
surface to have 
its Principal Configuration  
Structurally
Stable, with respect to $C^3 $ small deformations of its immersion into $\mathbb R^3$, are the following:

    \ref{th.1}.a) The umbilic points must all be 
    Darbouxian.

     \ref{th.1}.b) The periodic curvature lines (principal cycles) must be all
    hyperbolic.  That is,  their  Poincar\'e  Transformation,
    or First Return Map, must have a derivative different of  $1$.

     \ref{th.1}.c) It must not admit  connections, or self-connections,   of umbilical separatrices.

     \ref{th.1}.d) The limit sets of any non-periodic principal curve  must be an umbilic point or a principal cycle.
\end{theorem}
\begin{figure}[H]
\centerline{
\includegraphics[scale=0.55]{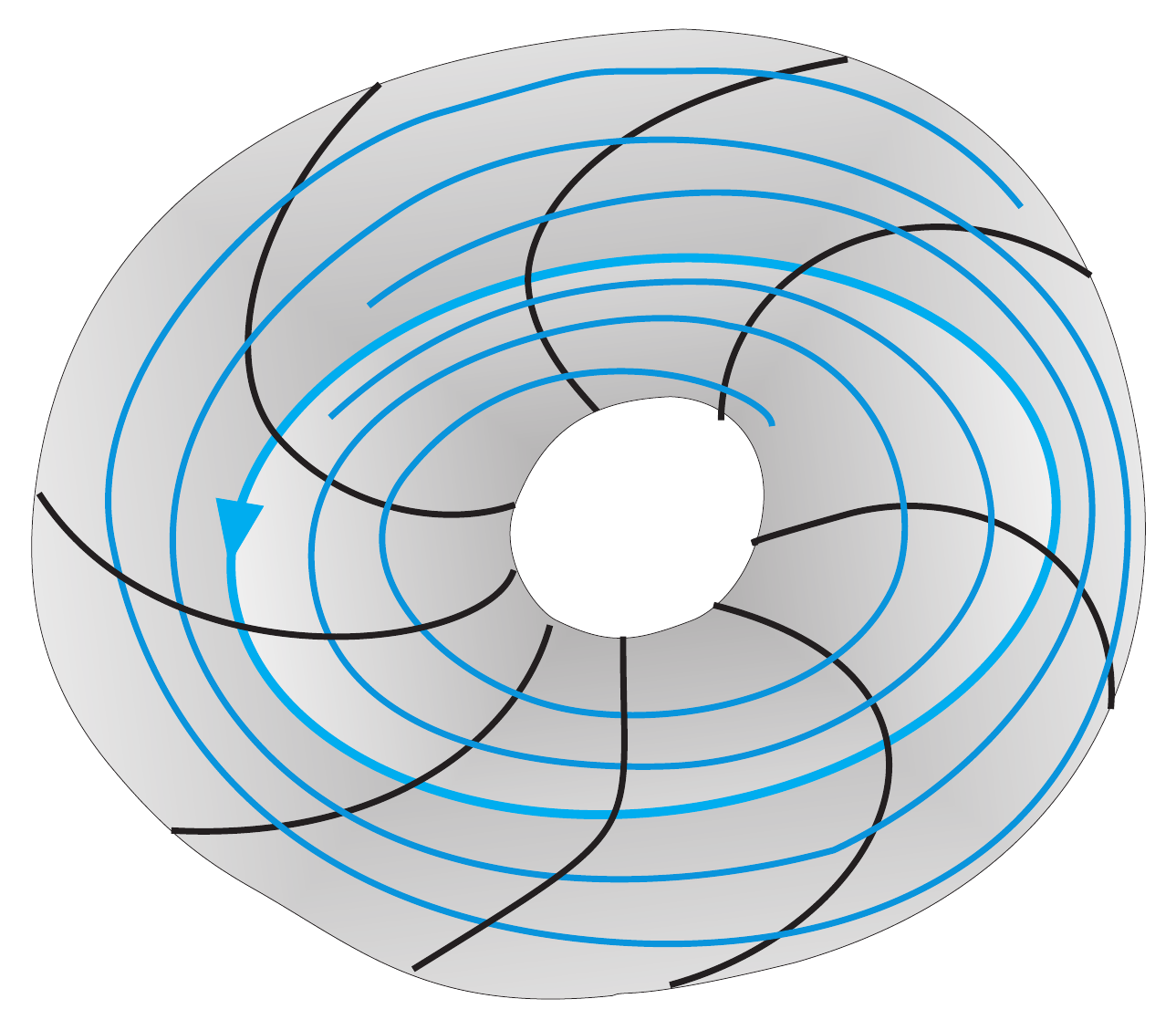}}
\caption{{\small  A periodic principal curvature  line, i.e. a principal cycle, which is hyperbolic. Illustration of  lines of curvature neighboring a  
   principal curve,  with first return map with  derivative
less than $1$. The arrows are placed to indicate a
local conventional orientation.  
The lines of curvature, 
generally,  are not globally orientable.  See, for example,
a neighborhood of the  
umbilic points in Fig. \ref{fig:5}.}} \label{fig:cycle}
\end{figure}
\begin{theorem} \label{th.2}
Within the space of immersions of  compact oriented surfaces,
endowed with the  $ C^ 3$ topology, the surfaces that satisfy the
conditions (\ref{th.1}.a) to (\ref{th.1}.d) form a set that is

     \ref{th.2}.a) open
     
      and
      
     \ref{th.2}.b) dense.

\end{theorem}

A rough scheme of demonstration for Theorems \ref{th.1} and \ref{th.2}, revealed 
 that part  (\ref{th.2}.a) of Theorem \ref{th.2}, as well as  the sufficiency of
conditions (\ref{th.1}.a) to (\ref{th.1}.d) for Theorem \ref{th.1}, seemed feasible.
However, 
serious 
difficulties  in establishing global aspects
of the necessity of the  conditions and, above all, the density, part
(\ref{th.2}.b) of Theorem \ref{th.2}, became  apparent.

At  
that
 time no 
 example  was  known of  a 
compact surface 
satisfying the conditions  (\ref{th.1}.a) to (\ref{th.1}.d).

The influence that  Peixoto's  work on the genericity of the 
{\it Structurally Stable  vector fields on two-dimensional manifolds,} 
Topology, 1962, 
was crucial  in the formulation of  the Fundamental Problem \ref{prob:2}.   This matter has been evoked  in the
 author's   essay:\\
 \noindent  $\color{blue}{\bullet}$ {{\it On a list of ordinary differential equations problems}, S\~ao Paulo Journal of Mathematical Sciences, 
 \url{https://doi.org/10.1007/s40863-018-0110-3}\label{lista}.\;  Zbl 1417.01036 \; Zbl 1028.34001. \;
MR3947401.  
 \noindent  $\color{blue}{\bullet}$  The    article 
{\tiny 
 \url{https://www.
 maths.ox.ac.uk/about-us/departmental-art/theory/differential-geometry},
 } 
 in an  Oxford University webpage,
 contains  a  very objective outline of Differential Geometry, with pertinent links.

\vskip 0.2cm
\begin{figure}[H]
	\begin{center}
\includegraphics[scale=0.37]{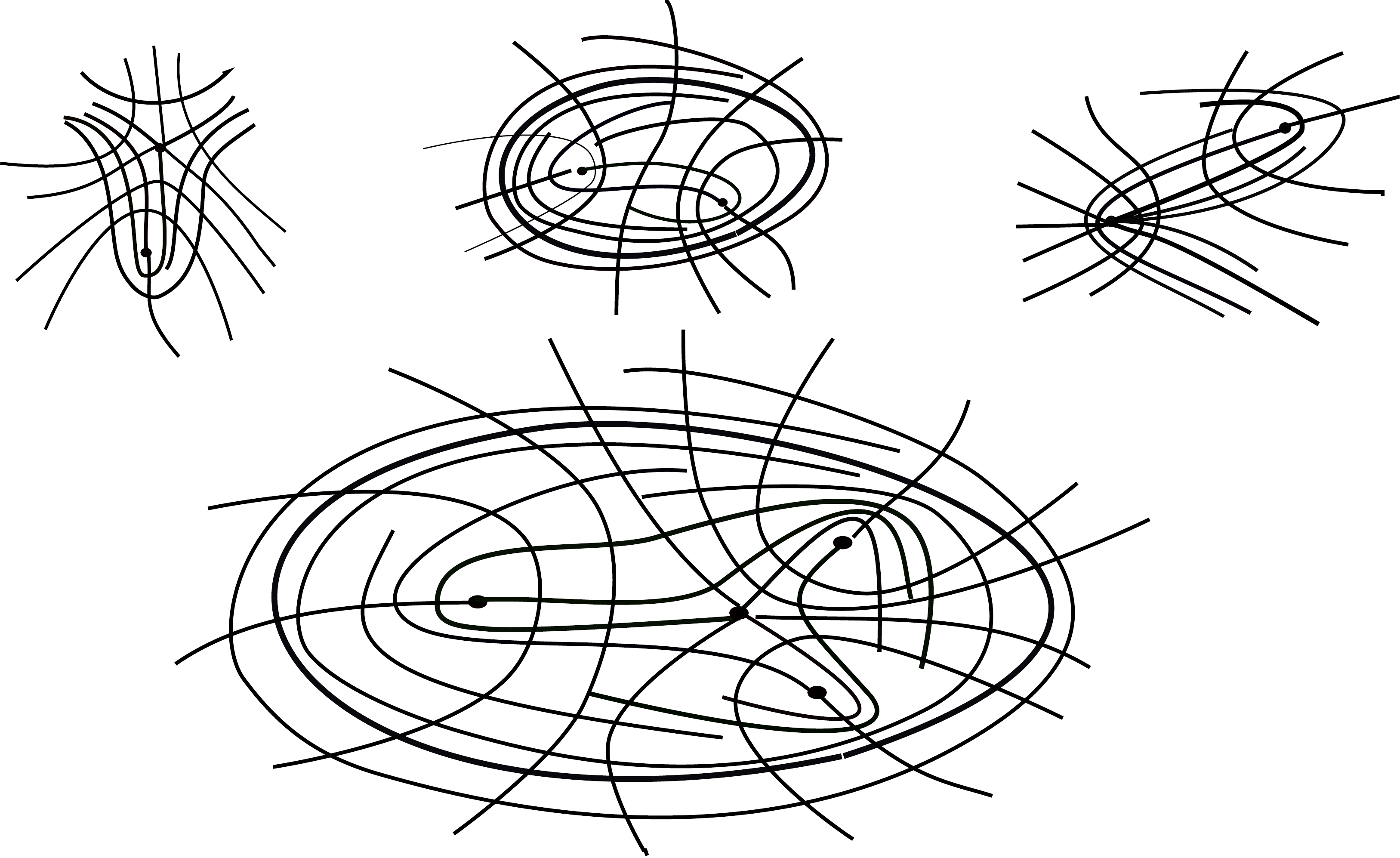}
 \caption{ Illustrations of  principally structurally stable patterns.  \label{fig:lcglobal} }
	\end{center}
\end{figure}

 \section{  Principal  Curvature 
 Configurations   as  a 
  Research  Area. 
 } \label{sc:3} 
A library search for
examples led nowhere.  
Not 
a single 
clue for  an  isolated principal  cycle, much
less hyperbolic, was encountered. 

The same for 
surfaces with 
recurrent principal curvature lines,  that is with those lines  that  violate  condition (\ref{th.1}d) in Theorem \ref{th.1},  
as it  is the case of integral curves of   non singular  vector fields with {\it irrational rotation number} on the 
torus.

The  Qualitative Theory of Differential Equations
 and   Dynamical Systems Theory,  founded by Poincar\'e, almost 100 years before,
had  not
penetrated this aspect of  Principal Configurations  in   Classical Differential Geometry, 
contrasting  with the progress that dynamic ideas had
in the study of geodesics.

A first round of calculations led  to a formula for the derivative 
 of the Poincar\'e First Return Map, in terms of the Mean Curvature
 $ H
= (k_1 + k_2) / 2 $, the Gaussian Curvature: $ K
= k_1 k_2 $ 
and the Christoffel Symbols. The elimination of these symbols  
was 
achieved 
later 
using   Codazzi  integrability conditions.

 The derivative, $T '$, of  the first return map, $ T $,  for  a  periodic curvature line
$\gamma $, after simplification achieved the 
following  expression: 
$$\log(T') = \pm {\frac{1}{2}}\int_{\gamma}   dH/
(H^2 -K)^{1/2}.$$

I 
accumulated a 
large  
collection 
of bibliographic 
references 
on 
principal curvature lines and umbilic points, 
none of them elucidating  Problem \ref{prob:2}. 
Among them  were  several  papers
relative to the   Ca\-ra\-th\'eo\-do\-ry Conjecture:

\begin{conjecture}\label{conj:cara}
Every $C^2$ convex and compact surface has at least two
umbilic points.
\end{conjecture} 

Years later I
found out 
 that this conjecture was considered
demonstrated for the case of analytic surfaces. 
 The conjecture is acknowledged to be  
open for $ C^\infty $   surfaces.

For the generic case it is trivial. In fact, having
the Darbouxian umbilic points indices  $ \pm 1/2 $, there must be at least
four of them  to be able to add up to give  the Euler-Poincar\'e  Characteristic of the
surface (in this case equal to $ 2 $).

Intrigued about the circumstances that could have led to 
the  formulation of this astonishing conjecture, with the aid of German colleagues,  I researched for 
its sources.     However, no  written record of this
formulation by the author was found. 
It is presumed that this was  
done orally
in seminars\footnote {The biography of C. Carath\'eodory by the historian Maria Georgiadou, Springer, 2004, 
suggests  to have  elucidated  the origin of this conjecture. 
However, the  reference  given  leads to a work on umbilic  points unrelated to the conjecture-problem \ref{conj:cara}.}.

I learned  that Gullstrand, the author for me initially
unknown, had been awarded   the  1911 Nobel Prize in Medicine, 
for his contribution to Ophthalmology. 
The Acta Mathematica 
paper  that I  
had browsed through 
was the mathematical part of the work 
that made him worthy of the award.

Discussions about aspects of my project  with  colleagues, more experienced than me,  in the fields of 
Dynamical Systems and Differential Geometry,   led to no mathematical reward. 
Unfortunately, Peixoto was traveling, outside Rio de Janeiro, during  that crucial semester.

I was lucky to meet  
the distinguished young  American Geometer Herbert Blaine Lawson, visiting IMPA.
He  was very receptive and  made 
the following  
thoughtful   
and  stimulating comment:
 
 ``{\it If what you propose works, 
 it will be opening a new research 
area.}"

 \section{ A  
 promenade with Monge, Darboux and Gullstrand. 
  } \label{sc:visits}
\subsection {S\~ao Paulo.} \label{ssc:sp}
In  March 1971,   I visited  the University of S\~ao Paulo.    I was  asked  by colleagues  to deliver  a
seminar lecture  presenting   the subject of  my research.  

This prompted me   
to organize the  pertinent material  I had gathered, 
 starting with the classical background: Euler, Monge, Dupin, then continuing with 
the geometric transformations that preserve the principal configurations: rigid motions, inversions and small  
parallel  displacements along the normals to the surface, 
etc. 
Concluding  with a discussion of the statements of Theorems \ref{th.1} and \ref{th.2}, the  Darbouxian Umbilics and, 
to include a  personal contribution,   the
formula for the derivative of the Poincar\'e return map  that  I had
found.  

Although  the 
seminar lecture on Principal Configurations did not happen 
because of administrative reasons, I could have  rewarding  discussions of  local  aspects of the project with Waldyr Oliva 
and Edgard Harle\footnote{\url{https://www.ime.usp.br/instituto/nossos-mestres/carlos-harle}} (1937 - 2020). 

In particular, Edgard Harle,  Geometer and fluent in  German, 
 helped me to 
 understand  
the article of 
  J. Fischer, Deutsche Math., 1935, written in the Gothic alphabet,  giving an  example of 
a principal 
configuration with  spiraling behavior  
around  an umbilic point. 
 
 The surfaces that verified the conditions (\ref{th.1}.b) and
(\ref{th.1}.d)   would be of this sort, around principal cycles.

Note that in the case of quadrics, surfaces of revolution
and all the other examples established in   Classical Differential Geometry, 
by the fact of deriving from a calculation, there is always
a first integral.


The visit to S\~ao Paulo and  the discussions, along the  preparation for a  seminar  presentation, helped me to reinforce the conviction that 
the subject was not devoid of interest.

\subsection {Salvador.} \label{ssc:ba}
 
The International Colloquium on Dynamical Systems of Salvador arrived,
in July  1971. There I met the  
French mathematician Ren\'e
Thom {1923 - 2002)  with whom I discussed  the expectations I had around The Fundamental Problem \ref{prob:2}.

It was an interesting dialog. 

For him, the umbilics
represented catastrophes, within the focal set: the envelope of 
the family of normal lines to the surface, that is  the caustic of the surface. 

When 
I asked something about the umbilical separatrices associated with the
points of Darboux, he   answered mentioning  properties  about the
``ridges" \, associated with the hyperbolic  and
elliptic,  (focal)  umbilics  of the generic surfaces. 

When I mentioned the nets of
lines of principal curvature on  surfaces, he answered me
with the focal set located in the ambient space.
 
 \subsection {Trieste and Paris.} \label{ssc:trsrt}
In July 1972, just before traveling to Trieste, to the International Center for Theoretical
Physics (ICTP), and to Paris,  for a visit to  the Institut de Hautes \'Etudes
Scientifiques (IHES), 
Carlos Guti\'errez,  a  doctoral
candidate at IMPA,  asked me for a suggestion of 
a thesis problem. 
 Without thinking twice, I took out of my briefcase
a copy of Darboux's article (which always accompanied me) and
I proposed: 

 \begin{problem}  \label{prob:gutierrez}
 Give a modern proof of
  this theorem, so that
mortals can understand it.\\
 Prove that 
conditions (\ref{th.1}.a) to (\ref{th.1}.d) in  the Fundamental Problem, Theorems \ref{th.1} and \ref{th.2},   are generic.  \newline
Pay special attention to 
the last condition  (\ref{th.1}.d).
 \end{problem}

On that trip I carried  with me a good part of the bibliographic material
relevant to Problem \ref{prob:gutierrez}, in  case I had to correspond with
Guti\'errez. 

To every  Swedish mathematician I met  in Trieste  I commented about Gullstrand paper. 
I inquired for the medical journal  in which he 
published his contribution to Ophthalmology. 
More than one colleague  promised to mail  me a copy,  which never arrived
 to my hands. 

By the end of November,  I received a letter from
Guti\'errez. 
He had decided to consider a more general form of
 Line Fields in two-dimensional manifolds. In that context,
he already had a coherent theory involving a weak form of $C^1$ density  for a class of Line Fields defined  by conditions resembling, topologically,  those formulated  in  (\ref{th.1}a)  to  (\ref{th.1}d), Theorem \ref{th.1}.

 His topologically  general theory, however, did not bring anything  enlightening on the genericity of the principal line fields and  their singularities (the umbilic points), as stated in the 
Fundamental Problem \ref{prob:2}.

I encouraged him to continue on the path he had set out,
which I also found interesting.

When we met again in January 1973, we quickly discussed the
difficulties encountered in attacking the  problem incorporating the  modification he
 suggested.  We
went on   considering  the  more abstract  approach, though quite  distant form Classical Differential Geometry.

In this direction he  converged to writing a  Doctoral Thesis,  IMPA 1974. 
In this endeavor, Guti\'errez became an exceptional expert in recurrent integral  curves of line fields. 
A fact  that will be crucial
for the unfolding, in  section \ref{sc:1981},  of  the story  recounted in this essay.  

 For  several years we did not speak again about  umbilic  points  and  principal  curvature configurations.

 \section{ Beaune   and  the  Hotel Monge  in Dijon.} \label{sc:beaune}

During the second semester of 1975 I visited Dijon for the first time.
After perambulating 
  several hotels, I 
ended lodged 
at the
small   Monge Hotel\footnote{Presently the  restaurant  ``Le Marrakech"  operates  in the building that hosted the hotel at 20, rue Monge.  See 
 \url{https://www.yelp.com/biz/le-marrakech-dijon-2}\label{monge_h_marrakech}  and    \url{https://www.yelp.com/biz/hotel-monge-dijon}.\label{monge_h_now}.}, located on the street with the same illustrious  name,
in the heart of the charming neighborhood called ``{\it Dijon Hist\'orique"}.
I stayed there for a week.

My host,   Robert Roussarie from the University of Dijon, took me to
visit  the famous vineyards of the Burgundy ``{\it C\^ote d'Or".}

We extended the visit to the historic village of Beaune. There, after a 
tour by impressive medieval buildings, I 
came upon, face to
face,  with an imposing statue of Gaspard Monge, 
{\small \url{https://commons.wikimedia.org/wiki/File:Beaune_-_Monument_de_Gaspard_Monge.jpg}
}, who was originally from that
town. 

The picture of the ellipsoid 
encountered  in 1970 in  Rio de Janeiro,  as 
evoked in section \ref{sc:1},    reappeared in my memory.
The cultural tour had 
excited my imagination. 
I was intrigued with the 
coincidence of my encounter with the statue and with the name of my 
hotel. 
The memory of the latent related problems \ref{prob:1}, \ref{prob:2} and  \ref{prob:gutierrez}, 
came back.

The  coincidence of historical, cultural and mathematical  circumstances intrigued me along the rest of my visit to Dijon, while being  lodged at the Hotel Monge.
  I profited from the occasion to browse Monge's book
\emph{Une application d'analyse \`a  la g\'eom\'etrie, 1795,}  at  the 
 University of Dijon Library, in the section of rare publications.     
 There he studies the general properties of  principal curvatures  and direction fields. 
 After perusing some of its sections, 
 I got convinced that as a  representative of the European Mathematics of the eighteenth century,
Monge would have had a hard time to 
understand  my gibberish and problematics, in \ref{prob:1}, \ref{prob:2} and  \ref{prob:gutierrez}, typical of an
anonymous tropical redoubt of the decade of the sixties, in the 20th century.

The link 
 {\small \url{https://books.google.com.br/books?id=aSEOAAAAQAAJ&redir_esc=y}}, leads pre\-sent\-ly to the above cited book of Monge.

 \section{
 Plastic  Transparency and Extrapolation  Exercises.} \label{sc:plastic}

Toward 1976, I had the help of mathematicians fluent in German, that translated for me substantial portions of Gullstrand's  paper  and  helped me  to  interpret  his  results. 
It was clear that the  author  was  strongly  interested in  the focal  normal set, specially in how its  ridge singular curves approached the umbilic points. 

However, he  did not address 
the analytic 
foundations to justify the principal configuration pictures of Darboux and, much less,
the other more degenerate ones that he also considered in his work. 
Subjects such as the uniqueness of the
umbilical separatrices, for example, were not touched at all.
 
 I  packed     the
 bibliographic material  pertinent  to the Fundamental Problem \ref{prob:2},  Theorems  \ref{th.1}   and   \ref{th.2},   and 
 also my 
 handwritten 
 notes,  
  in a plastic  bag, 
 which  I  carried  up and down with me. 
For long periods, when my briefcase was too heavy, I used to  leave   the bag 
resting at home.  Months later I would  take it along with me  again.

The transparency of the plastic prompted me to locate and have   the Fundamental Problem 
present. 
The  transportation ceremony  ritual gave me a  strange sensation of 
possession   and of being
working on it,  secretly.

With uncertain periodicity, however,   I  engaged 
in longer 
intuitive  exercises, 
without directly  attacking  its
fundamental mathematical difficulties. \newline 
I  particularly  enjoyed extrapolating  conditions  \ref{th.1}.a to 
\ref{th.1}.d,  weakening and adapting  them  to higher order ones,  predicting  
 the  generic bifurcations of principal curvature configurations  on surfaces evolving subject to one   parameter spatial deformations.
 This amounted to an additional  {\it 
 stage in the sequence of qualitative jumps} discussed in sections \ref{sc:1} and \ref{sc:2}.  
 
A bolder, higher order,   extrapolation  to which the  previous predictions pointed to  were  
the  principal configurations of  $3-$dimensional manifolds   immersed into $\R^4$.  
 
 In section \ref{sc:books} will be given references to works achieved later in the  lines of research outlined above. 

 \section{ Collaboration in Research Teams.}\label{sc:teams}
At the beginning of 1980 we received a determination from the CNPq Central Administration: 
Researchers should be organized in
teams, as in soccer, and present joint projects. 
This
was a novelty.  We were accustomed to receiving  
requests  to fill out new cadastre  forms every time
there was a change in the administrative cadres.

This naive administrative measure ended up being beneficial for
me, for 
Carlos Guti\'errez and also for the Fundamental Problem on
Principal Curvature Configurations.  

I engaged  Guti\'errez to join forces with me to work in the
Issue. 

The progress was surprisingly fast, taking into account that
each of us worked simultaneously on 
individual projects.
The two of us,  and the Problem,  had gained maturity.\newline
We updated the formula for the derivative, $T '$, of the first return transformation, $ T $, of a periodic curvature  line 
$ \gamma $.

In terms of the mean $ H $ and Gaussian $ K $ curvatures, 
the beautiful following expression, mentioned also in section \ref{sc:3},    holds:

$$\log(T') = \pm {\frac{1}{2}}\int_{\gamma}   dH/
(H^2 -K)^{1/2}.$$

In fact, chronologically,  it was obtained after  its equivalent version 
$$ \log (T ') = \pm \int_{\gamma}  dk_2 /(k_2-k_1), $$
 in terms of  the principal curvatures
 $ k_1 <k_2 $,  since $ H = (k_1 + k_2) / 2 $ and $ K = k_1k_2$, 

\noindent from 
which we derived   a perturbation  method to 
hyperbolize  periodic curvature lines. 

Using a  method of resolution (``blowing up") of singularities, very close to the 
traditional one, used to reduce the study of 
complicated singular points to simpler ones, we fully justified,   for  the case of surfaces of  class $ C^4$,  the
Darbouxian principal configuration pictures  in Fig. \ref{fig:5}. 

This was a novelty in relation to the analytic case considered
by the French Geometer.

On the recurrent principal  lines of curvature, Carlos Guti\'errez, 
formulated  the following diagnosis:\newline

{\it It will be possible to get an approximation in class $C^2$ that eliminates them.
 The increasing from  class  $ C^2$ to $C^3$, however,  will be very difficult.} \newline
 
We produced  an  example of a Toroidal Surface  in $\mathbb R^3$ without umbilic points
and with dense recurrent curvature lines. Its immersion, however, was quite distant from the standard Torus of revolution.

The pieces of the puzzle seemed to fit.  
We had a
sustainable version!

The time had come to make a broad communication of the
results. The 1981 International Dynamical Systems Symposium inaugurating the
new installations of IMPA, constituted an 
auspicious occasion. 
 
\section{The  
first  international presentation  on Principal Curvature Configurations.
} \label{sc:1981}

 On the eve of the 
 lecture,  I stayed until
later than usual  discussing with Guti\'errez what the presentation would be like. Some 
delicate points 
emerged.

 It  became explicit  
 that we did not have 
 examples of principal 
recurrences other than that in the aforementioned, free of umbilics,  Torus.
In particular,
we did not know 
immersions of  spheres, necessarily carrying umbilic points, exhibiting principal recurrences.
Furthermore, at  that time the Torus  example seemed  too technical
to be  quickly explained  geometrically during the
lecture.

We were in the uncomfortable position of being able to eliminate all
the recurrent  curvature lines, replacing them after
an approximation $C^2$-small, by periodic lines of curvature or, in
some cases, by  connections of umbilical separatrices. 
Specifically we knew only the example of the Torus.

What if someone, interested in specific  examples, raised the question? 
Worst still, what if there were no more examples of recurrences than those on the Torus?

In that case, our results would be considerably  weakened.

The search for an example of principal recurrences on a surface of
genus  zero lasted a long while. 

It was then that  Monge's Ellipsoidal   reappeared;   this time  it was  more flexible than ever. 
It allowed more daring deformations and contortions, including
non-analytic ones.
So, bending more here,
rotating there,  we 
stumbled over  an example.

The 
lecture
 was ready!  It would be the first one on the 
next morning.

I started projecting  
Monge's Ellipsoid (Fig. \ref {fig:1}) and
reviewed how to   explain  its principal configuration, using the Theorem of
Dupin (Fig. \ref{fig:4}).
I  commented on  the absence of examples 
of Global  Principal Configurations, 
distinct from 
that one.

After the preparatory setting concerning the space of Immersions of oriented two-manifolds into $\R^3$,  I proposed the problem of recognizing globally the immersion with Structurally Stable  principal configurations and  establishing their  genericity.
Then came the statement of  our result,
insisting on the limitation in class $C^2 $ for the density 
approach. 
I ended up summoning the listeners to raise this
class to  $ C^3 $, thus solving the only problem that remained
open in relation to the initial program, as stated in Theorems \ref{th.1} and \ref{th.2}. 

\begin{remark}  \label{prob:2to3}
This problem remains open until now.
\end{remark} 

Reaching the end of the lecture, none of the experts in recurrences on  Dynamical Systems 
 in the audience  asked for specific  examples. This was  disappointing. 

Fortunately,  Dan Henry\footnote{\url{https://www.ime.usp.br/map/dhenry/danhenry/texto01.htm}} (1945 - 2002), an expert  in
Differential Equations in Infinite Dimensions, formulated the
expected question:

``{\it I do not know why you make the hypothesis (\ref{th.1}.d). I do not know any
surface that does not satisfy  it.}"

I drew an ellipsoid of revolution (Fig. \ref{fig:2}, left). I deformed it slightly
so that,  in both polar  caps,  it was like one of Monge, with
their three axes distinct, while around the equator it continued
being 
one of revolution.

Then, I rotated around the major axis only the
upper hemisphere of the surface. 
Because of its 
equatorial rotational symmetry, this was a family of $C^\infty$
surfaces $ E_\theta $, depending  on the parameter $\theta$,  
designating  the angle of rotation.
\begin{figure}[H]
\begin{center}
\includegraphics[scale=0.7]{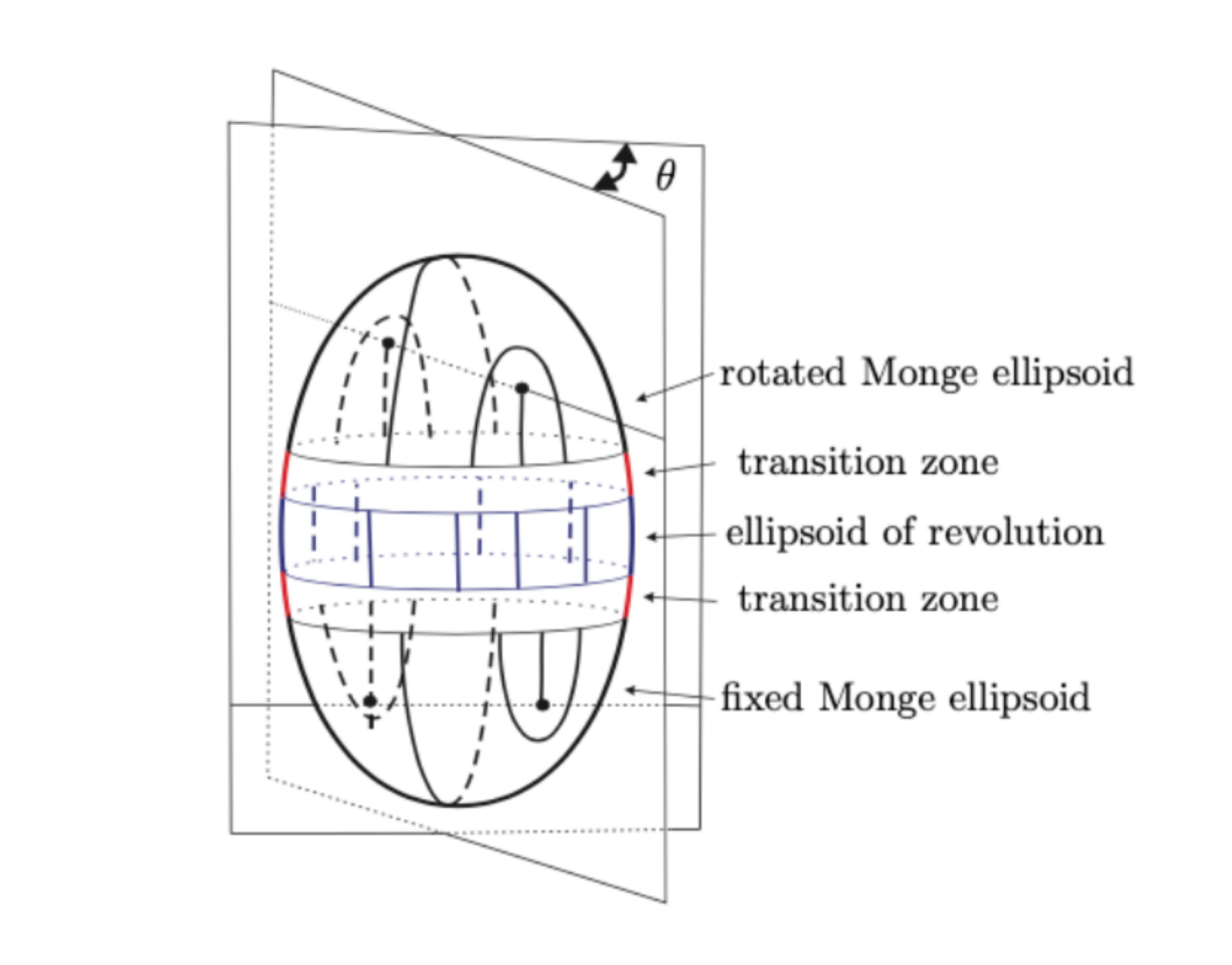}
\end{center}
\caption{ {\small Ellipsoidal Surface $E_\theta$ with oscillatory recurrent principal curvature  lines.} \label{fig:9} }
\end{figure}

It is evident that the curvature lines of $E_\theta $
define in their second return to the equatorial circle a rotation
of angle $2\theta $.
 So, for  angles incommensurable with
respect  to $ 2\pi $, the lines of curvature become all dense
in $ E_\theta $. See illustration in Fig. \ref{fig:9}.

The  written  presentation of the work was divided  in two parts. 

The first one, establishing that conditions (\ref{th.1}.a)  to 
(\ref{th.1}.d) define a  $C^3$  open  set consisting in Structural Stability immersed surfaces, was published in {\it
Asterisque}, Vol. 98-99, 1982. 

The second part, which
demonstrates that any compact and oriented surface can be
arbitrarily $ C^2$ approximated  by  one that verifies
the four above mentioned  conditions, appeared in {\it Springer
Lecture Notes in Mathematics}, Vol. 1007, 1983.

\subsection{Record of the  first mathematical encounter between Principal Curvature  Geometry and Structural Stability.} \label{ss:encounter} 
 These papers are pointed out as the  documentation pertinent to  the first  encounter between  the line of thought disclosed 
   from the works of   Monge, 1796, Dupin, 1815,  and Darboux,  1896,   with that transpiring from the achievements  of Poincar\'e, 1881, 
   Andronov - Pontrjagin, 1937, and Peixoto, 1962.
 
\vspace{-0cm}
\begin{figure}[H]
\begin{center}
\includegraphics[scale=0.7]{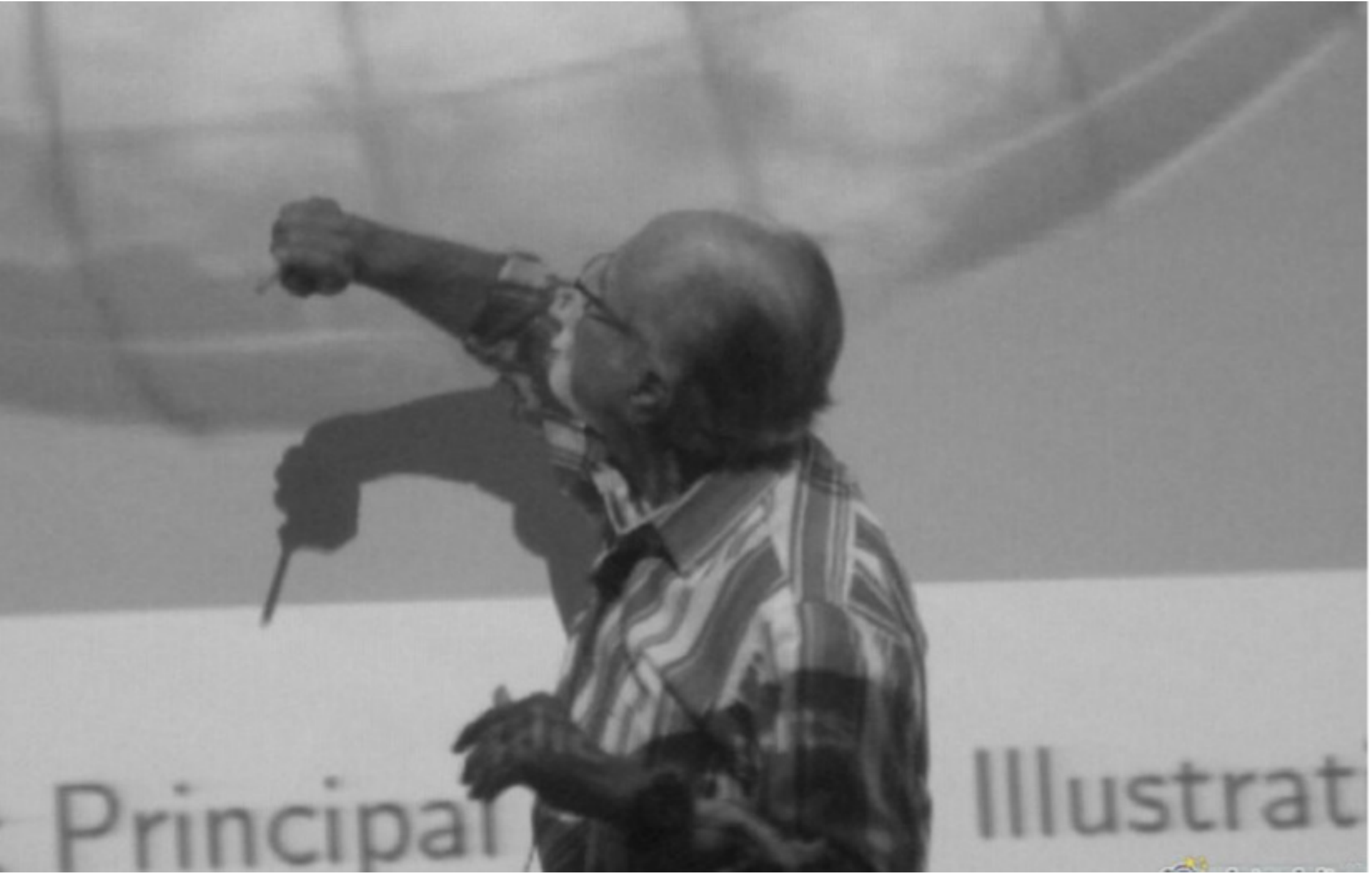}
\end{center}
\caption{{\small No photo of the lecture at IMPA, 1981, was found.  Here is one from a lecture on Principal Configurations  delivered in Montevideo, 2012. 
Photo by A. Chenciner. The ellipsoid picture was borrowed from   E. Ghys site.} \label{fig:10}}   
 \end{figure}

\section{ A Herald of the Qualitative Theory of Differential Equations.
} \label{sc:herald}

In the first half of 1990, on the occasion of the  ``{\it Ann\'ee Sp\'eciale
de  Syst\`emes Dynamiques",}   sponsored by the CNRS (National Research Council) of France, 
I delivered  a short course   at the
University of Dijon. The subject was {\it  Umbilic Points and Principal 
Curvature Lines.}
 
During  that visit I had access to the original work of
Monge, entitled  {\it Sur les lignes de courbure de la surface de
l'Ellipsoide}, published in {\it Journal de l'Ecole
Polytechnique., II cah, 1796}.

The contemplation of Monge's original pictures, 
supplemented 
with  some reflections and additional readings,  led me to
make explicit 
an intuitive observation 
about the history of 
mathematical ideas, which, in an embryonic form, was already present in
my initial motivation, after the first contact with the Ellipsoid and the
subsequent formulation of the Fundamental Problem \ref{prob:1}.

There is no bibliographic record that Euler, responsible for the
conceptualization of the principal  curvature  line fields,
would have integrated them, visualizing the Principal  Configuration.

Monge was the first mathematician to recognize the importance of this
structure, providing the first non-trivial example, for
global integration of the differential equations of the lines
of curvature (that is, of the line fields $ L_1 $ and $ L_2 $) in the
case of the Ellipsoid. 

Let's inquire:
{\it What would have led him to establish a result of this
nature, which in our days fits perfectly within the Qualitative Theory of 
Ordinary Differential  Equations,
 considering that this happened almost one hundred years before
 Poincar\'e founded it and defined its goals?}
Monge's motivation 
derived from a complex interaction of
aesthetic and practical considerations, and also the  explicit intention 
to apply the results of his mathematical research. 

At
Monge's time was more tenuous the artificial separation between
Mathematics and Applications. 
The lines of curvature were discovered
studying the problem of transportation  of 
debris 
in the
construction of embankments and fortifications,  called
the  {\it ``deblais et remblais"}   problem. 
His ellipsoidal picture was proposed 
in the architectural project 
for  the construction of the vault,  over an elliptical terrain,   of 
Legislative Assembly Building  of  the French Revolution Government: the lines of
curvature would be the guiding curves  for the  placement of the stone bricks, the
umbilics would serve as supporting  points 
for hanging
 the sources of illumination, under one of which would be
located  the rostrum for   the speakers.
The architectural  project was never
put into effect.  
\begin{remark}\label{rm:herald}
 If  Poincar\'e, for the scope and depth of its contribution,
is recognized as the founder, Monge has the merit of being regarded as a 
precursor --a herald-- of the Qualitative Theory of Differential Equations
and  of Foliations with Singularities Theory.
\end{remark}

However, it should be emphasized  that the novelty of the 
Qualitative Theory of Poincar\'e lies in the methods developed  for 
the study of the phase portrait of general non-integrable differential equations, 
that he established for 
polynomial  equations, in  the generic case.
Keeping in mind that this was a laboratory model for the far  reaching and more complex problems 
of Celestial Mechanics that he was investigating. \newline

The connection
between the works of Monge and Poin\-ca\-r\'e, outlined above, 
was overlooked in the accomplished 
historical work of Ren\'e Tat\'on  ``{\it L'Oeuvre Scientifique de
Monge,}"   Presses Univ.  de France, 1951.

\section{ Books, later developments  and updated  
references.}\label{sc:books}

 \subsection{A Short Course in 1991.} \label{ssc:IMPA_91}
The  project  of writing an
expository  book based on the lectures delivered     at  the  1990 short course in  Dijon, was  formulated on that occasion. 
See  the beginning of  section \ref{sc:herald}.

 In 1991, I engaged  
Guti\'errez  to  collaborate with me in 
writing  a small book of didactic vocation, 
proposed for the 
18th Brazilian Mathematics Colloquium. 

Besides   the results
on Structural Stability and Approximation  established  in  the original papers,
the lecture notes 
included  more details concerning  the
theoretical foundations and the  motivation for the study of Principal Curvature Configurations.
 The examples of principal 
recurrences were improved and reformulated in more 
conceptual terms, especially the one devoted to   the  immersed Toroidal Surface\footnote{This example of 1991 has only one of the principal foliations with dense curves.   In  {\sc R. Garcia} and {\sc J. Sotomayor,}   Tori embedded in $\mathbb R^3$ with dense principal lines. Bull. Sci. Math.,  {\bf  133:4 }  (2009),  348-354,  was given an example in which  both principal  foliations have its lines  dense.}.
The  bibliographic  references were also updated. 

  \noindent  $\color{blue}{\bullet}$ {\sc  C. Gutierrez  } and {\sc J. Sotomayor,}
  {Lines of Curvature and Umbilical Points on Surfaces}, $18^{\hbox{th}}$  Brazilian  Math.
Colloquium, Rio de Janeiro,  IMPA, (1991). \newline 
  Reprinted as
 { Structurally Stable  Configurations of Lines of Curvature and Umbilic
 Points on Surfaces, Lima, Monografias del IMCA}, (1998). 
MR2007065.

\subsection{A book  in 2009.} \label{ssc:IMPA_09}
A book of   broader scope than that of 1991,  with more topics on  the rich interaction of Classical Differential Geometry and Differential Equations,
 was published in 2009.
 \newline 
  \noindent  $\color{blue}{\bullet}${\sc R.~Garcia} and {\sc J.~Sotomayor}, \newblock{\it Differential Equations of Classical Geometry,
  	a Quali\-ta\-tive Theory},  
	Publica\c c\~oes Matem\'{a}ticas, 27$^{\text \small{o}}$ Col\'oquio Brasileiro de Matem\'atica, IMPA, (2009).   Zbl 1180.53002. 
MR2532372.

 \subsection{Further developments on surfaces in $\R^3$.} \label{ssc:M2R3}
Without claiming to be complete, some
additional works that are 
pertinent  to the  
 present 
essay are listed below in chronological order.

\vspace{0.15cm}
 \noindent$\color{blue}{\bullet}$ {\sc J. W.  Bruce} and {\sc D. L.  Fidal,} { \em On binary differential equations and umbilic
points,} Proc. Royal Soc. Edinburgh  {\bf 111A}, 1989, 147-168. Zbl 0685.34004.

  \noindent$\color{blue}{\bullet}$ {\sc C. Gutierrez } and {\sc  J.  Sotomayor,}
 {\em  Principal lines on surfaces immersed with constant mean
 curvature.}
  Trans. Amer. Math. Soc. {\bf 293}, 1986, no. 2, 751 - 766.  MR816323. Zbl 0598.53007.
  
    \noindent$\color{blue}{\bullet}$ {\sc  R. Garcia}  and {\sc  J. Sotomayor,} 
{\em  Lines of curvature near singular points of implicit surfaces}. Bull. Sci. Math. 117 (1993), no. 3, 313 - 331.  
MR1228948.

 \noindent$\color{blue}{\bullet}$ {\sc  R. Garcia}  and {\sc  J. Sotomayor,}
 {\em Lines of curvature near hyperbolic principal cycles}. Dynamical systems (Santiago, 1990), 255 - 262,
  Pitman Res. Notes Math. Ser., 285, Longman Sci. Tech., Harlow, 1993. MR1213951.
  
\noindent$\color{blue}{\bullet}$ {\sc C. Gutierrez } and {\sc  J.  Sotomayor,}
 {\em Periodic lines of curvature bifurcating from Darbouxian umbilical connections}. 
 Bifurcations of planar vector fields (Luminy, 1989), 196 - 229, Lecture Notes in Math., 1455, Springer, Berlin, 1990.   MR1094381.

  \noindent$\color{blue}{\bullet}$  {\sc  R. Garcia}  and {\sc  J. Sotomayor,}
 {\em Lines of Curvature on Algebraic Surfaces},
 Bull. Sciences  Math. {\bf  120}, (1996),   367-395. MR1411546.
 
 \noindent$\color{blue}{\bullet}$  {\sc  J. Sotomayor,}
 {\em Lines of curvature and an integral form of Mainardi-Codazzi equations.} An. Acad. Brasil. Ci\^enc. 68 (1996), no. 2, 133 - 137.
 MR1751266.
 
  \noindent$\color{blue}{\bullet}$   {\sc  T. Maekawa, F. E. Wolter} and {\sc  N. M. Patrikalakis,}  { \em Umbilics and lines of curvature for shape interrogation,} 
     Comput. Aid. Geometr. Des. {\bf 13}  (1996),  133 --161.  Zbl 0875.68858. 

   \noindent$\color{blue}{\bullet}$ {\sc  R. Garcia   } and {\sc J.  Sotomayor,} {\em Structural stability of parabolic
    points and periodic asymptotic lines},
 Matem\'atica Contempor\^anea, {  \bf 12},   (1997),  83-102.  
MR163442. 

  \noindent$\color{blue}{\bullet}$ {\sc C. Gutierrez} and {\sc  J.  Sotomayor,}
 {\em Lines of Curvature,  Umbilical Points and Ca\-ra\-th\'eo\-do\-ry Conjecture},
  Resenhas IME-USP,  {\bf 03}, 1998, 291-322.  
MR1633013. 

  \noindent$\color{blue}{\bullet}$ {\sc  R. Garcia}, {\sc   C. Gutierrez } and {\sc  J. Sotomayor,}
     {\em Lines of principal curvature around umbilics and Whitney umbrellas.}  Tohoku Math. J. (2)  {\bf 52:2}  (2000),   163-172. MR1756092. 

 \noindent$\color{blue}{\bullet}$ {\sc R. Garcia } and {\sc  C.  Gutierrez, } { \em   Ovaloids of $\mathbb R\sp 3$ and their
umbilics: a differential equation approach.}   J. Differential Equations {\bf 168}, 2000, no. 1, 200--211. MR1801351.

   \noindent$\color{blue}{\bullet}$ {\sc R. Garcia } and {\sc  J. Sotomayor,} {\em Structurally stable configurations of
   lines of mean curvature and  umbilic points on surfaces immersed in ${\mathbb R}^3$,}
     Publ. Matem\'atiques.  { \bf 45:2},  (2001),    431-466. MR1876916.  Zbl 0875.68858.
   
    \noindent$\color{blue}{\bullet}$ {\sc R. Garcia } and {\sc  J. Sotomayor,} { \em  Umbilic and tangential
singularities on configurations of principal curvature lines.} An.
Acad. Brasil. Ci\^enc. {\bf 74},  2002, no. 1, 1--17.  
MR1882514. 

   \noindent$\color{blue}{\bullet}$  {\sc R. Garcia } and  {\sc  J. Sotomayor,}  { \em Lines of Geometric Mean Curvature
   on surfaces immersed in ${\mathbb R}^3$,}
    Annales de la Facult\'e des Sciences de Toulouse, {\bf 11}, (2002),  377-401. 
MR2015760. 
 
 \noindent$\color{blue}{\bullet}$ {\sc R. Garcia } and {\sc  J. Sotomayor,} { \em Lines of Harmonic Mean
	Curvature on surfaces immersed in ${\mathbb R}^3$,}
Bull. Bras. Math. Soc.,  {\bf 34}, (2003), 303-331. 
MR1992644. 

 \noindent$\color{blue}{\bullet}$ {\sc R. Garcia } and {\sc  J. Sotomayor,} { \em Lines of  Mean
   Curvature on surfaces immersed in ${\mathbb R}^3$,}
   Qualit. Theory of Dyn. Syst. {\bf 5}, 2004,  137-183. MR2129722.

   \noindent$\color{blue}{\bullet}$ {\sc  R. Garcia}, {\sc   C. Gutierrez } and {\sc  J. Sotomayor,} { \em  Bifurcations of Umbilic Points and Related
  	Principal Cycles.}  Journ. Dyn. and Diff. Eq. {\bf 16, 2}, (2004), 321-346. MR2129722. 

 \noindent$\color{blue}{\bullet}$ {\sc R. Garcia } and {\sc  J. Sotomayor,} { \em  On the patterns of principal curvature lines around a curve of umbilic points.}
  An. Acad. Brasil. Ci\^enc. 77 (2005), no. 1, 13 - 24.  MR2114929 

  \noindent$\color{blue}{\bullet}$  {\sc R. Garcia} and {\sc  J. Sotomayor}, 
      { \em Lines of principal curvature near singular end points of surfaces in $\mathbb R^3$.}  
      Singularity theory and its applications, 437-462, Adv. Stud. Pure Math., {\bf 43}, Math. Soc. Japan, Tokyo, (2006).
      MR2325150. 

  \noindent$\color{blue}{\bullet}$  {\sc R. Garcia} and {\sc J. Sotomayor,} 
   { \em Tori embedded in $\mathbb R^3$ with dense principal lines.} Bull. Sci. Math.,  {\bf  133:4 }  (2009),  348-354. 
MR2503006. 

  \noindent$\color{blue}{\bullet}$  {\sc R.  Garcia, R. Langevin} and {\sc P  Walczak},   {\em Foliations making a constant angle with principal directions on ellipsoids.}
Annales Polonici Mathematici,  {\bf 113}  (2015), 165-173.  
MR3312099.

\noindent$\color{blue}{\bullet}$  {\sc R.~Garcia} and {\sc J.~Sotomayor}, { \em Historical Comments on Monge's Ellipsoid and the Configurations of Lines of Curvature on Surfaces,}   Antiquitates Mathematicae,   {\bf 10(1)}  (2016),  348-354. Zbl 1426.53008.  \; MR3613151.

\noindent$\color{blue}{\bullet}$ {\sc  V. V. Ivanov,}   {\em An analytic conjecture of
Carath\'eodory. } Siberian Math. J. {\bf 43}, no. 2,  2002,
251--322.  MR1902826.

 \noindent$\color{blue}{\bullet}$ {\sc B. Smyth} and {\sc  F. Xavier, }
{\em  A sharp geometric estimate for the index of an umbilic on a
smooth surface.}  Bull. London Math. Soc. {\bf 24}, 1992, no. 2,
176--180.  MR1148679.

\noindent$\color{blue}{\bullet}$ {\sc B. Smyth} and {\sc  F. Xavier, }  {\em  Eigenvalue
estimates and the index of Hessian fields.}  Bull. London Math.
Soc. {\bf 33}, 2001, no. 1, 109--112.  MR1798583.

\noindent$\color{blue}{\bullet}$ {\sc B. Smyth,}  {\em The nature of elliptic sectors in
the principal foliations of surface theory.}  EQUADIFF 2003, 957--959, World Sci. Publ., Hackensack, NJ, 2005. 

\subsection{Principal Configurations of Hypersurfaces in $\R^4$.} \label{ssc:M3R4}
An important 
achievement that 
must  be mentioned is the extension of
the generic properties of the lines of curvature associated to  
three-dimensional manifolds  immersed in $\mathbb R^4 $,
first obtained by Ronaldo Garc\'ia in his Thesis (IMPA, 1989). See partial list below. \newline
\noindent$\color{blue}{\bullet}$   {\sc R. Garcia } {\it
Principal  Curvature Lines near Partially Umbilic Points in
hypersurfaces immersed in $ {\mathbb R}^ 4 $,  Comp. and Appl. Math.,
20}, 121 - 148, 2001.

 In this direction,  the following more recent developments should be added: \newline
 \noindent$\color{blue}{\bullet}$  {\sc  D. Lopes, J. Sotomayor } and R. {\sc Garcia}, \newblock{\it 
 	Umbilic and Partially umbilic singularities   of
 	Ellipsoids  of  $\mathbb R^4$}. Bulletin of the Brazilian Mathematical Society, {\bf 45}, (2014),  453-483. 
	 
 \noindent$\color{blue}{\bullet}$   {\sc  D. Lopes, J. Sotomayor } and R. {\sc Garcia}, \newblock{\it Partially umbilic singularities   of
 	hypersurfaces of  $\mathbb R^4$}. Bulletin de 
 Sciences Mathematiques, {\bf 139}, (2015),  421- 472. 
 
  \noindent$\color{blue}{\bullet}$   {\sc R. Garcia } and {\sc  J. Sotomayor,}  Lines of curvature on quadric hypersurfaces of $\mathbb R^4$, Lobachevskii Journal of Mathematics {\bf 37}, (2016), 288-306. 

 \subsection{Axial Configurations  on surfaces in $\R^4$.} \label{ssc:M2R4} 

 The name {\it axial} refers to the association with the axes of the {\it normal ellipse of curvature} 
 defined at the regular points  of surfaces mapped into  $\R^4$.

   \noindent$\color{blue}{\bullet}$ {\sc  R. Garcia   } and {\sc J.  Sotomayor,}
 {Lines of axial curvature on surfaces immersed in $\mathbb R^4$}, Differential Geom. Appl.
 {  \bf 12}, (2000),  253--269.   MR1764332. \; Zbl 0992.53010.

\noindent$\color{blue}{\bullet}$ {\sc L. F.  Mello,}
 {\em Mean directionally curved lines on surfaces immersed in ${\Bbb R}\sp 4 ,$}
 Publ. Mat. { \bf  47}, (2003), 415--440.
 
\noindent$\color{blue}{\bullet}$  {\sc R. Garcia, L. Mello} and {\sc J. Sotomayor},  { \em Principal mean curvature foliations on surfaces immersed in
       $\mathbb R^4$,}   EQUADIFF (2003), 939-950, World Sci. Publ., Hackensack, NJ, 2005.  \; Zbl 1116.57024.

 \noindent$\color{blue}{\bullet}$   {\sc R. Garcia}  and {\sc   J. Sotomayor},  { \em Lines of axial curvature at critical points  on surfaces mapped into  
    $\mathbb R^4$,}  Sao Paulo J Math Sci   {\bf 6},  (2012), no. 2,  277- 300.

  \noindent$\color{blue}{\bullet}$   {\sc R. Garcia,  J. Sotomayor}  and {\sc F. Spindola},  { \em Axiumbilic  Singular Points on Surfaces Immersed in $\mathbb R^4$ and their Generic Bifurcations,}   Journal of Singularities  {\bf 10}, (2014), 124-146.

\subsection{Comments on some open problems raised  in the works listed in this essay}\label{ssc:open_problems}

\subsubsection{Increase  in the density class from $ C^2$  to $ C^3$,    as inquired  in  Theorem \ref{th.2} } 

\subsubsection{Recurrence on cubic deformations of Monge's  ellipsoid. 
From exercise 3.6.3, p. 83,  of  the 2009, Garcia and Sotomayor book.} 

 For $\rho$ small,  consider the cubic surface $S_\rho = f_\rho ^{-1}(0)$                                     
$$ f_\rho (x,y,z) = \frac{x^2}{a^2} +  \frac{y^2}{b^2} +z^2 +\rho xyz - 1;  \; \;  a>0, \;  b>0,  (a - 1)(b - 1)(a - b) \neq 0.$$
 From simulations of the possible global behaviors of principal foliations of $S_\rho,$
 formulate a conjecture about the possibility of dense principal lines on algebraic surfaces of spherical type.
 
 The surface $S_\rho,$ for $ \rho \neq 0$ is the simplest algebraic one that exhibits the twisting effect depicted in the smooth example in Fig. \ref{fig:9}.  This effect, determined by the cubic term  $\rho xyz$, was first studied in Garcia and Sotomayor, Bull. Sci. Math. 1993, for the case of the  cubic surface 
 $$  \frac{x^2}{a^2} +  \frac{y^2}{b^2} - z^2 +\rho xyz =0;  \; \;  a>0, \;  b>0,  (a - 1)(b - 1)(a - b) \neq 0, $$
 locally conoidal  at  its critical point $(0,0,0)$.
It also appears in  the study of the general case of principal configurations on  algebraic surfaces. See  Garcia and Sotomayor, Bull. Sci. Math. 1996.

No mathematical proof of the existence of recurrences  on $S_\rho,$ for $ \rho \neq 0$  has been published yet.

 \subsection{Two pertinent presentations.} \label{ssc:youtube}
    Two lectures 
     delivered by the  author,    overlapping with  the  subject
     of  this   essay   have been  posted in:
    \newline           
\noindent$\color{blue}{\bullet}$   \url{https://youtu.be/JX2pHiCvaxw}\label{soto_peixoto_90},  in Workshop International  de Sistemas Din\^amicos -
 $90$ Mauricio Peixoto,  2011, and in  \newline
\noindent  $\color{blue}{\bullet}$ 
 \url{https://youtu.be/IscUm7UHv50}\label{soto_gutierrez_impa_18},  A lecture evocative of Carlos Guti\'errez, 
2018.  
 \begin{figure}[H]
\begin{center}
\includegraphics[scale=0.8]{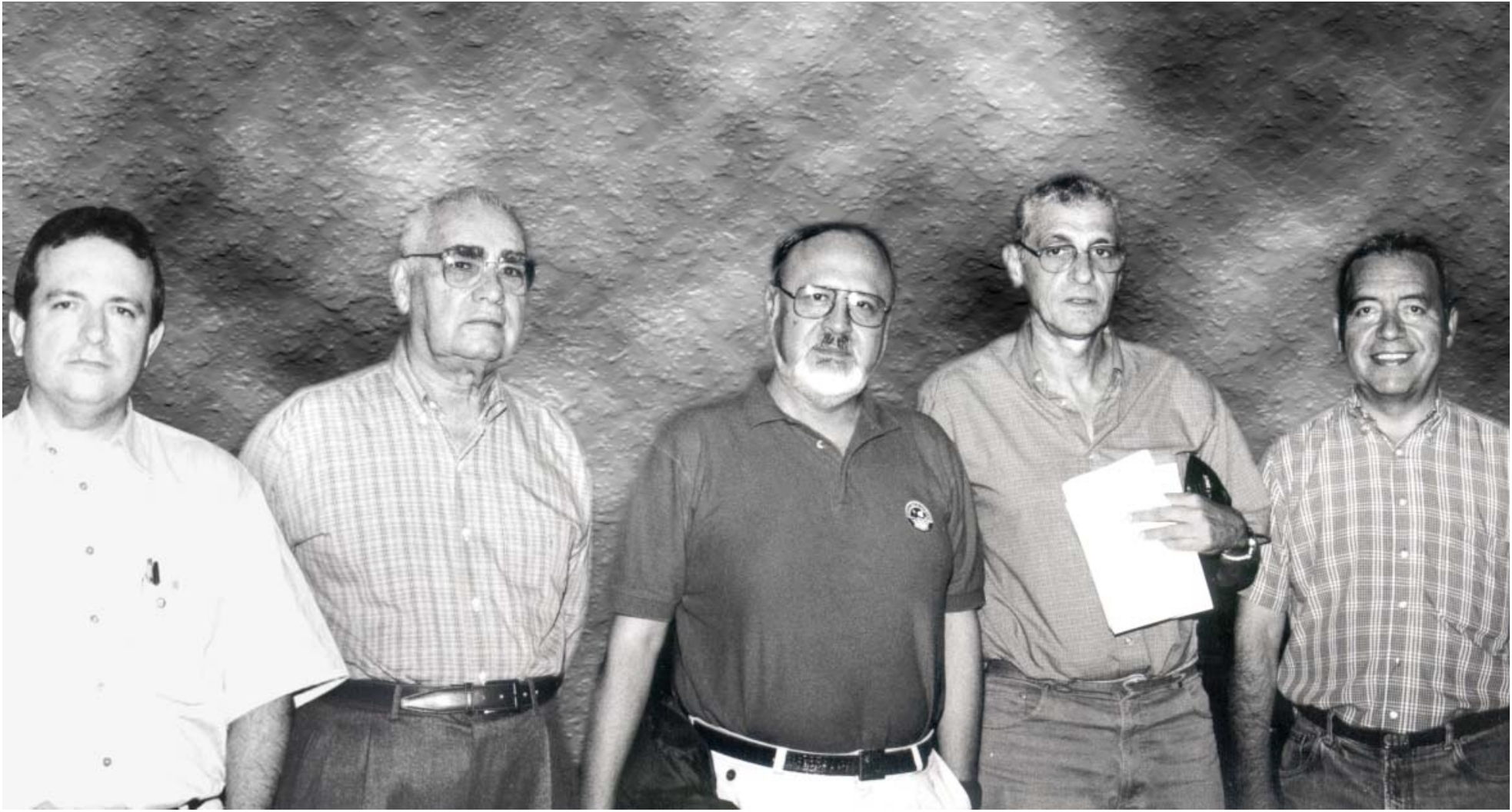}
\end{center}
\vspace {-0cm}
\caption{{\small  At UNICAMP, 
 2000, three mathematical generations: right to left,  C. T.  Guti\'errez, M. A. Teixeira, J.  Sotomayor, M. M. Peixoto and R. A. Garcia}. \label{fig:11}}   
\end{figure} 

\subsection{Closing words evoking       
Poincar\'e:}
\begin{verse} \label{avenir}
{\it If we wish to foresee the future of Mathematics, \\
the right approach is to study 
 its history and\  present condition.\\ 
For us mathematicians,  is it not this procedure, to some extent, routine?\\
We are accustomed to extrapolation, which is a method \\
of deducing the future from the past and the present;\\
and,  since we are well aware of its limitations, we run no risk  \\
of deluding ourselves as to the scope of the results it provides.} 
\end{verse}
H. Poincar\'e, in The Future of Mathematics  (\emph{L'Avenir des Math\'ematiques}), 
speech read by G. Darboux at the  ICM,   Rome, 1908.

 \subsection{On the genesis of this essay.} \label{ssc:elipsoide}
 This work is an English version,  which includes  considerable  revision, upgrading and adaptation,   based on \\
 \noindent$\color{blue}{\bullet}$ \emph{O Elips\'{o}ide de Monge,}  Portuguese, Mat. Univ., 15, 1993,
and  its Spanish  translation in Materials Matem\`atiques,
2007:  \url{http://mat.uab.cat/matmat/PDFv2007/v2007n01.pdf}.

\vspace{0.15cm}

 \noindent{\bf Acknowledgements.\,}  The author is  grateful to  M. O. Sotomayor, C. P. Moromisato,  F. E. Wolter and 
  L. F. Mello    for  helpful  style suggestions on a previous version   and to R. A. Garcia for his mathematical comments and substantial aid in 
  the production of the pictures in this  work.
 
 
\textcolor{white}{  
}
 
 \noindent \author{\noindent Jorge Sotomayor\\ Instituto de Matem\'atica e Estat\'{\i}stica,\\  
Universidade  de S\~ao Paulo,\\
 Rua do Mat\~ao  1010. \\  
Cidade Univerit\'aria, CEP 05508-090,\\
S\~ao Paulo, S. P., Brazil}\; \\

\end{document}